\pdfoutput=1
\RequirePackage{ifpdf}
\ifpdf 
\documentclass[pdftex]{sigma}
\else
\documentclass{sigma}
\fi

\numberwithin{equation}{section}

\newtheorem{Theorem}{Theorem}[section]
\newtheorem{Corollary}[Theorem]{Corollary}

 { \theoremstyle{definition}

 }

\newcommand{\ct}{\operatorname{CT}}

\begin{document}

\allowdisplaybreaks

\renewcommand{\thefootnote}{}

\newcommand{\arXivNumber}{2212.02355}

\renewcommand{\PaperNumber}{059}

\FirstPageHeading

\ShortArticleName{A New (But Very Nearly Old) Proof of the Rogers--Ramanujan Identities}

\ArticleName{A New (But Very Nearly Old) Proof\\
 of the Rogers--Ramanujan Identities\footnote{This paper is a~contribution to the Special Issue on Basic Hypergeometric Series Associated with Root Systems and Applications in honor of Stephen C.~Milne's 75th birthday. The~full collection is available at \href{https://www.emis.de/journals/SIGMA/Milne.html}{https://www.emis.de/journals/SIGMA/Milne.html}}}

\Author{Hjalmar ROSENGREN}

\AuthorNameForHeading{H.~Rosengren}

\Address{Department of Mathematical Sciences, Chalmers University of Technology\\ and University of Gothenburg, SE-412~96 G\"oteborg, Sweden}
\Email{\href{mailto:hjalmar@chalmers.se}{hjalmar@chalmers.se}}
\URLaddress{\url{http://www.math.chalmers.se/~hjalmar/}}

\ArticleDates{Received March 25, 2024, in final form June 27, 2024; Published online July 02, 2024}

\Abstract{We present a new proof of the Rogers--Ramanujan identities. Surprisingly, all its ingredients are available already in Rogers seminal paper from 1894, where he gave a~considerably more complicated proof.}

\Keywords{Rogers--Ramanujan identities; constant term identities}

\Classification{11P84; 33D15}

\begin{flushright}
\emph{Dedicated to Stephen Milne on his 75th birthday}
\end{flushright}

\renewcommand{\thefootnote}{\arabic{footnote}}
\setcounter{footnote}{0}

\section{Introduction}
The fascinating story of the Rogers--Ramanujan identities
\begin{subequations}\label{rr}
\begin{align}\sum_{k=0}^\infty\frac{q^{k^2}}{(1-q)\bigl(1-q^2\bigr)\dotsm(1-q^k)}&=\prod_{j=0}^\infty\frac 1{\bigl(1-q^{5j+1}\bigr)\bigl(1-q^{5j+4}\bigr)},\label{rra}\\
\sum_{k=0}^\infty\frac{q^{k(k+1)}}{(1-q)(1-q^2)\dotsm\bigl(1-q^k\bigr)}&=\prod_{j=0}^\infty\frac 1{\bigl(1-q^{5j+2}\bigr)\bigl(1-q^{5j+3}\bigr)}\label{rrb}
\end{align}
\end{subequations}
has been told many times. They were first proved by Rogers in a brilliant but forgotten paper from 1894 and rediscovered by Ramanujan some time before 1913.
Hardy wrote about the seven proofs known in the 1930s: ``None of them is both simple and straight-forward, and probably it would be unreasonable to hope for a proof which is'' \cite{h}.
Several dozen more proofs have appeared since then, using diverse methods.
We refer to \cite{a89} for a survey of the proofs available in the 1980s and \cite{akp,bp,c,f,gis,p,st,t} for examples of more recent work.
The purpose of the present note is to give a new proof that, although far from the shortest possible, we find relatively ``simple and straight-forward''.
Almost unbelievably, this proof was very nearly found by Rogers. He knew all its main ingredients and could easily have given it by adding
one or two sentences to \cite{r2}.
Among the many gems buried in Rogers' paper are the identities \cite[Section~6, formulas~(5)--(8)]{r2}
\begin{subequations}\label{ghr}
\begin{alignat}{3}
&G(q)+G(-q)=2F(q)G\bigl(q^{16}\bigr),\qquad&&
G(q)-G(-q)=2qF(q)H\bigl(-q^4\bigr),&\\
&H(q)+H(-q)=2F(q)G\bigl(-q^4\bigr),\qquad&&
H(q)-H(-q)=2q^3F(q)H\bigl(q^{16}\bigr),&
\end{alignat}
\end{subequations}
where $G(q)$ and $H(q)$ denote the sums in \eqref{rra} and \eqref{rrb}, respectively, and
\begin{equation*}
F(q)=\prod_{j=0}^\infty\frac{1-q^{8j+8}}{1-q^{2j+2}}.
\end{equation*}
 Below, we will give both a version of Rogers' proof of \eqref{ghr} and another proof based on work of Andrews \cite{a86}.
After deriving \eqref{ghr}, Rogers writes somewhat dismissively
\begin{quote}These four identities \dots\ are sufficiently remarkable in themselves to call for mention at this point, although they may all be derived from the $\theta$-function values of the series $\phi(q)$, $\psi(q)$ [our $G(q)$ and $H(q)$] obtained in the last section. [\dots] It is not, however, in these identities that the special interest in the series $\phi(q)$ and $\psi(q)$ lies.
\end{quote}
At this point Rogers has already proved the Rogers--Ramanujan identities. He realizes that
 \eqref{ghr} are quite easy to prove for the product sides of \eqref{rr} and hence perhaps less interesting than they may seem.

Our observation is simply that the identities \eqref{ghr} uniquely determine $G$ and $H$ up to normalization. In fact, if they are written as
\begin{subequations}\label{abf}
\begin{gather}
G(q)=F(q)\bigl(G\bigl(q^{16}\bigr)+qH\bigl(-q^4\bigr)\bigr), \\
H(q)=F(q)\bigl(G\bigl(-q^4\bigr)+q^3H\bigl(q^{16}\bigr)\bigr),
\end{gather}
\end{subequations}
 they can be used to recursively compute $G$ and $H$ as power series, starting from the initial values $G(0)=H(0)=1$. Since Rogers has proved that \eqref{ghr} hold both for the sum sides and the product sides of \eqref{rr}, this gives an independent proof of the Rogers--Ramanujan identities, which is
considerably easier than the one given earlier in the same paper.

We will present the resulting proof in some detail, but we stress that apart from the preceding paragraph everything can be found in
 \cite{a86,r2}.
I am very much indebted to George Andrews for pointing this out to me.

In Section~\ref{ois}, we explain how some further known Rogers--Ramanujan-type identities either appear as byproducts of the proof or can be obtained by slight variations.
We only consider identities that are very close to \eqref{rr}, in that the product sides are essentially the same. It would be interesting to investigate whether similar proofs can be given for other results. As a first step,
in \cite{ro} we show that the method can be used to prove
Rogers' and Slater's mod $7$ identities.

\section{A proof of the Rogers--Ramanujan identities}\label{ps}

We will use the standard notation
\begin{equation*}
(a;q)_k=\prod_{j=0}^{k-1}\bigl(1-aq^j\bigr),\qquad k\in\mathbb Z_{\geq 0}\cup\{\infty\}, \end{equation*}
which is extended to negative indices by $(a;q)_{-k}=1/\bigl(aq^{-k};q\bigr)_k$. More generally,
\begin{equation*}(a_1,\dots,a_m;q)_k=(a_1;q)_k\dotsm (a_m;q)_k. \end{equation*}
We will write
\begin{equation*}R(t;q)=\sum_{k=0}^\infty\frac{q^{k^2}t^k}{(q;q)_k} \end{equation*}
and $G(q)=R(1;q)$, $H(q)=R(q;q)$.
The Rogers--Ramanujan identities \eqref{rr} then take the compact form
\[
G(q)=\frac 1{\bigl(q,q^4;q^5\bigr)_\infty},\qquad H(q)=\frac 1{\bigl(q^2,q^3;q^5\bigr)_\infty}.
\]
Finally, we introduce the multiplicative theta function $\theta(z;q)=(q,z,q/z;q)_\infty$.

The only results that we will need are Euler's identity
\begin{equation}\label{esa}\frac 1{(z;q)_\infty}=\sum_{k=0}^\infty\frac{z^k}{(q;q)_k},\qquad |z|<1,\end{equation}
 Jacobi's identity
\begin{gather}\theta(z;q)=\sum_{k=-\infty}^\infty(-1)^kq^{k(k-1)/2}z^k \label{jti}
\end{gather}
and
\begin{gather}
\frac{\theta(1/z;q)}{(tz;q)_\infty}=(t;q)_\infty\sum_{k=-\infty}^{\infty}\frac{(-1)^kq^{k(k-1)/2}}{(t;q)_kz^k}, \qquad 0<|z|<\frac 1{|t|}.\label{rps}
\end{gather}
To some readers these identities may look just as complicated as \eqref{rr}, but they are much easier to prove. For instance, they can all be obtained as special cases of Ramanujan's ${}_1\psi_1$-summation, a simple proof of which was given by Ismail \cite{i}, see also \cite[Appendix~C]{a86}.

Following Andrews \cite{a86}, we will represent the sums in \eqref{rr} as constant terms. In general, if~${f(z)=\sum_{k=0}^\infty a_k z^k}$
 (all series in this paper can be considered as formal or analytic according to taste)
and we let~$\ct$ denote the constant term of a Laurent series in $z$, then by \eqref{jti}, we~have%
\begin{align}
\ct \theta(1/z;Q)f\bigl(xz^N\bigr)&=\ct \sum_{l=-\infty}^\infty (-1)^lQ^{l(l-1)/2}z^{-l}\sum_{k=0}^\infty a_k x^kz^{Nk}\nonumber\\
&=\sum_{k=0}^\infty
(-1)^{Nk}Q^{Nk(Nk-1)/2}a_k x^k.\label{fexp}
\end{align}
We apply this with $f(z)$ as Euler's series \eqref{esa}. In order to get a factor $q^{k^2}$, we must choose $Q=q^{2/N^2}$. We will only need the cases $N=1$ and $N=2$, which give the constant term representations
\begin{equation*}R(t;q)
=\ct\frac{\theta\bigl(1/z;q^{2}\bigr)}{ (-qtz;q)_\infty}
=\ct\frac{\theta\bigl(1/z;q^{1/2}\bigr)}{ \bigl(q^{1/2}tz^2;q\bigr)_\infty}.
\end{equation*}

Next, we write
 \begin{gather*}
 \frac{\theta\bigl(1/z;q^{2}\bigr)}{ (-qtz;q)_\infty}= \frac{\theta\bigl(1/z;q^{2}\bigr)}{\bigl(-qtz;q^2\bigr)_\infty}\cdot\frac1{\bigl(-q^2tz;q^2\bigr)_\infty}=\frac{\theta(1/z;q^{2})_\infty}{\bigl(-q^2tz;q^2\bigr)_\infty}\cdot\frac1{\bigl(-qtz;q^2\bigr)_\infty},\\
\frac{\theta\bigl(1/z;q^{1/2}\bigr)}{ \bigl(q^{1/2}tz^2;q\bigr)_\infty}= \frac{\theta\bigl(1/z;q^{1/2}\bigr)}{ \bigl(q^{1/4}t^{1/2}z;q^{1/2}\bigr)_\infty}\cdot \frac 1{ \bigl(- q^{1/4}t^{1/2}z;q^{1/2}\bigr)_\infty},
\end{gather*}
and expand the resulting factors using \eqref{esa} and \eqref{rps}. This gives the quadratic transformation formulas
\begin{subequations}\label{rtf}
\begin{align}
R(t;q)&=\bigl(-qt;q^2\bigr)_\infty\sum_{k=0}^\infty\frac{q^{k(k+1)}t^k}{\bigl(q^2,-qt;q^2\bigr)_k}\\
&=\bigl(-q^2t;q^2\bigr)_\infty\sum_{k=0}^\infty\frac{q^{k^2}t^k}{\bigl(q^2,-q^2t;q^2\bigr)_k}\label{rtfb}\\
&=\bigl(q^{1/4}t^{1/2};q^{1/2}\bigr)_\infty\sum_{k=0}^\infty\frac{q^{k^2/4}t^{k/2}}{\bigl(q^{1/2},q^{1/4}t^{1/2};q^{1/2}\bigr)_k}.\label{rtfc}
\end{align}
\end{subequations}
These identities were found by Rogers \cite[Section~6, formulas~(2)--(4)]{r2}. Except for his way of writing, the original proof is not difficult; we will review it in Section~\ref{rpss}.
 Watson independently proved \eqref{rtfb} and \eqref{rtfc} \cite{w}.

For $t=1$ and $t=q$, the equations \eqref{rtf} simplify to
\begin{subequations}\label{rm20}
\begin{align}
G(q)&=\bigl(-q;q^2\bigr)_\infty\sum_{k=0}^\infty\frac{q^{k(k+1)}}{(-q;-q)_{2k}}\\
& =\bigl(-q^2;q^2\bigr)_\infty\sum_{k=0}^\infty\frac{q^{k^2}}{\bigl(q^4;q^4\bigr)_k} \\
&=\bigl(q^{1/4};q^{1/2}\bigr)_\infty\sum_{k=0}^\infty\frac{q^{{k^2}/4}}{\bigl(q^{1/4};q^{1/4}\bigr)_{2k}},\\
H(q)&=\bigl(-q^2;q^2\bigr)_\infty\sum_{k=0}^\infty\frac{q^{k(k+2)}}{\bigl(q^4;q^4\bigr)_k}\\
& =\bigl(-q;q^2\bigr)_\infty\sum_{k=0}^\infty\frac{q^{k(k+1)}}{(-q;-q)_{2k+1}} \\
&=\bigl(q^{1/4};q^{1/2}\bigr)_\infty\sum_{k=0}^\infty\frac{q^{{k(k+2)}/4}}{\bigl(q^{1/4};q^{1/4}\bigr)_{2k+1}}.
\end{align}
 \end{subequations}

Let us introduce the auxiliary series
\begin{equation}\label{ab}
A(q)=\sum_{k=0}^\infty\frac{q^{k^2}}{\bigl(q^4;q^4\bigr)_k},\qquad B(q)=\sum_{k=0}^\infty\frac{q^{k(k+2)}}{\bigl(q^4;q^4\bigr)_k}. \end{equation}
Then, \eqref{rm20} can be written
\begin{subequations}\label{ghij}
\begin{gather}
 G\bigl(-q^4\bigr)=\bigl(q^4;q^8\bigr)_\infty\frac{B(q)+B(-q)}{2},\label{ga}\\
G(q)=\bigl(-q^2;q^2\bigr)_\infty A(q),\\
 G\bigl(q^{16}\bigr)=\bigl(q^4;q^8\bigr)_\infty\frac{A(q)+A(-q)}{2}, \\
H(q)=\bigl(-q^2;q^2\bigr)_\infty B(q),\label{hb}\\
 H\bigl(-q^4\bigr)=\bigl(q^4;q^8\bigr)_\infty\frac{A(q)-A(-q)}{2q},\\
H\bigl(q^{16}\bigr)=\bigl(q^4;q^8\bigr)_\infty\frac{B(q)-B(-q)}{2q^3}.
\end{gather}
\end{subequations}
We may use \eqref{ga} and \eqref{hb} to eliminate
$A$ and $B$ from the remaining equations. Writing also
\begin{equation*}
F(q)=\frac{\bigl(-q^2;q^2\bigr)_\infty}{\bigl(q^4;q^8\bigr)_\infty}=\frac{\bigl(q^8;q^8\bigr)_\infty}{\bigl(q^2;q^2\bigr)_\infty},
\end{equation*}
we arrive at \eqref{ghr} or, equivalently, \eqref{abf}.

As was noted in the introduction, the identities \eqref{abf}
can be considered as recursions for computing the Taylor series of $G$ and $H$, starting from $G(0)=H(0)=1$; see \eqref{rec} below. Hence, if we can verify that
 \eqref{ghr} also holds for the product sides
\begin{equation*}
\tilde G(q)=\frac 1{\bigl(q,q^4;q^5\bigr)_\infty},\qquad
\tilde H(q)=\frac 1{\bigl(q^2,q^3;q^5\bigr)_\infty},\end{equation*}
then the Rogers--Ramanujan identities $G=\tilde G$, $H=\tilde H$ follow.
As Rogers observed, this reduces to classical theta function identities.
We first write
\begin{align*}
\tilde G(q)&=\frac 1{\bigl(q,q^4,q^6,q^9;q^{10}\bigr)_\infty}=\frac{\bigl(-q,-q^9;q^{10}\bigr)_\infty}{\bigl(q^4,q^6;q^{10}\bigr)_\infty\bigl(q^2,q^{18};q^{20}\bigr)_\infty}=\frac{\bigl(q^8,q^{12};q^{20}\bigr)_\infty \theta\bigl(-q;q^{10}\bigr)}{(q^2;q^{2})_\infty}
 \end{align*}
 and similarly
 \begin{equation*}
 \tilde H(q)=\frac{\bigl(q^4,q^{16};q^{20}\bigr)_\infty \theta\bigl(-q^3;q^{10}\bigr)}{\bigl(q^2;q^{2}\bigr)_\infty}. \end{equation*}
Thus, we are reduced to the problem of computing the even and odd part of a theta function.
In general, Jacobi's triple product identity \eqref{jti} implies the classical identities
\begin{align*}
&\frac{\theta(z;q)+\theta(-z;q)}2=\sum_{k=-\infty}^\infty q^{k(2k-1)}z^{2k}=\theta\bigl(-q z^2;q^4\bigr), \\
&\frac{\theta(z;q)-\theta(-z;q)}2=-\sum_{k=-\infty}^\infty q^{k(2k+1)}z^{2k+1}=-z\theta\bigl(-q^3z^2;q^4\bigr).
\end{align*}
It follows that, for instance,
\begin{equation*}\frac{\tilde G(q)+\tilde G(-q)}2=\frac{\bigl(q^8,q^{12};q^{20}\bigr)_\infty\bigl(-q^{12},-q^{28},q^{40};q^{40}\bigr)_\infty}{\bigl(q^2;q^{2}\bigr)_\infty}.\end{equation*}
Using that
\begin{equation*}
\bigl(-q^{12},-q^{28};q^{40}\bigr)_\infty =\frac {\bigl(q^{24},q^{56};q^{80}\bigr)_\infty}{\bigl(q^{12},q^{28};q^{40}\bigr)_\infty}
=\frac{\bigl(q^{16},q^{24};q^{40}\bigr)_\infty}{\bigl(q^{12},q^{28};q^{40}\bigr)_\infty\bigl(q^{16},q^{64};q^{80}\bigr)},
\end{equation*}
the expression above is easily reduced to
\begin{equation*}\frac{\tilde G(q)+\tilde G(-q)}2 =\frac{\bigl(q^8;q^8\bigr)_\infty}{\bigl(q^2;q^2\bigr)_\infty\bigl(q^{16},q^{64};q^{80}\bigr)_\infty}=\frac{\bigl(q^8;q^8\bigr)_\infty}{\bigl(q^2;q^2\bigr)_\infty} \tilde G\bigl(q^{16}\bigr).
\end{equation*}
The remaining parts of \eqref{ghr}, with $G$ and $H$ replaced by $\tilde G$ and $\tilde H$, follow similarly.
Alternatively, the identities \eqref{abf} can be obtained directly as special cases of
\begin{equation*}
\theta(z;q)=\theta\bigl(-qz^2;q^4\bigr)-z \theta\bigl(-q^3z^2;q^4\bigr).
\end{equation*}

Finally, we make some remarks on the
 recursions for the Taylor coefficients of $G$ and $H$ mentioned above.
Substituting
\begin{equation*}F(q)=\sum_{k=0}^\infty f_kq^{2k},\qquad G(q)=\sum_{k=0}^\infty g_kq^k,\qquad H(q)=\sum_{k=0}^\infty h_kq^k\end{equation*}
into \eqref{abf} gives
\begin{subequations}\label{rec}
\begin{alignat}{3}
& g_{2k}=\sum_{j=0}^{\lfloor k/8\rfloor}f_{k-8j}g_j,\qquad&&
g_{2k+1}=\sum_{j=0}^{\lfloor k/2\rfloor}(-1)^jf_{k-2j}h_j,&\\
& h_{2k}=\sum_{j=0}^{\lfloor k/2\rfloor}(-1)^jf_{k-2j}g_j,\qquad&&
h_{2k+1}=\sum_{j=0}^{\lfloor (k-1)/8\rfloor}f_{k-8j-1}g_j.&
\end{alignat}
\end{subequations}
It is interesting to note that all the coefficients involved here have partition-theoretic interpretations. This is straight-forward for $f_n$, and very well-known for $g_n$ and $h_n$.
This leads to the following result, which can be viewed as a combinatorial reformulation of~\eqref{abf}.

\begin{Corollary}
Let $f_n$ be the number of partitions of $n$ into parts not divisible by $4$.
Let $g_n$ be the number of partitions of $n$ into parts congruent to $\pm 1$ $\operatorname{mod}$ $5$ or, equivalently, the number of partitions
of $n$ without repeated or consecutive parts. Let $h_n$ be
the number of partitions of $n$ into parts congruent to $\pm 2$ $\operatorname{mod}$ $5$ or, equivalently, the number of partitions
of $n$ without repeated or consecutive parts and not containing $1$ as a part. Then, the recursions \eqref{rec} hold.
\end{Corollary}

\section{Some additional results}\label{ois}

Having proved the Rogers--Ramanujan identities, we can sum the six series in \eqref{rm20}. This was of course clear to Rogers, and leads to the identities
given in Slater's list \cite{msz,si,s} as number 99, 20, 98, 16, 94 and 96, respectively (as is noted in \cite{si}, 98 is equivalent to 79, and 94 to 17).

Above, we applied \eqref{fexp} to Euler's series \eqref{esa}. If we instead use the companion series
\[
 f(z)=(-z;q)_\infty=\sum_{k=0}^\infty\frac{q^{k(k-1)/2}z^k}{(q;q)_k},
 \]
we must take $Q=q^{1/N^2}$ to get a factor $q^{k^2}$. We will only consider the case $N=1$, which leads~to
\begin{gather}\label{rct}R(t;q)=\ct \theta(1/z;q) (qtz;q)_\infty.\end{gather}
Writing
\begin{align*}
\theta(1/z;q)(qtz;q)_\infty&=\frac{\theta(1/z;q)}{(-qtz;q)_\infty }\bigl(q^2t^2z^2;q^2\bigr)_\infty\\
&
= \frac{\theta(1/z;q)}{\bigl(q^{1/2}tz;q\bigr)_\infty}\bigl(q^{1/2}tz;q^{1/2}\bigr)_\infty=\frac{\theta(1/z;q)}{\bigl(q^{3/2}tz;q\bigr)_\infty}\bigl(qtz;q^{1/2}\bigr)_\infty \end{align*}
and expanding each factor gives the transformations
\begin{subequations}\label{rx}
\begin{align}
\label{ra}R(t;q)&=(-qt;q)_\infty\sum_{k=0}^\infty\frac{(-1)^kq^{3k^2}t^{2k}}{\bigl(q^2;q^2\bigr)_k(-qt;q)_{2k}}\\
\label{rb}&=\bigl(q^{1/2}t;q\bigr)_\infty\sum_{k=0}^\infty\frac{q^{k(3k-1)/4}t^k}{\bigl(q^{1/2};q^{1/2}\bigr)_k\bigl(q^{1/2}t;q\bigr)_k}\\
\label{rc}&=\bigl(q^{3/2}t;q\bigr)_\infty\sum_{k=0}^\infty\frac{q^{k(3k+1)/4}t^k}{\bigl(q^{1/2};q^{1/2}\bigr)_k\bigl(q^{3/2}t;q\bigr)_k}.
\end{align}
\end{subequations}
The identity \eqref{ra} is the limit case $B\rightarrow\infty$ of \cite[equation~(7.9)]{gs} and \eqref{rb} is
\cite[equation~(7.17)]{gs}. We have not found \eqref{rc} in the literature but it can be deduced from \eqref{rb} by a simple contiguity argument.
Namely, splitting the terms in \eqref{rc} using
\begin{equation*}
\bigl(q^{3/2}t;q\bigr)_\infty\frac{t^k}{\bigl(q^{3/2}t;q\bigr)_k}= \bigl(q^{1/2}t;q\bigr)_\infty\left(\frac{t^k}{\bigl(q^{1/2}t;q\bigr)_k}+\frac{q^{k+\frac 12}t^{k+1}}{\bigl(q^{1/2}t;q\bigr)_{k+1}}\right)
\end{equation*}
 and then replacing $k$ by $k-1$ in the second sum, one finds that \eqref{rc} equals \eqref{rb}.

In the cases $t=1$ and $t=q$, \eqref{rx} can be written as
\begin{subequations}
\begin{align}
\label{raa}
G(q)&=(-q;q)_\infty\sum_{k=0}^\infty\frac{(-1)^kq^{3k^2}}{(q^4;q^4)_k(-q;q^2)_{k}}\\
& =\bigl(q^{1/2};q\bigr)_\infty\sum_{k=0}^\infty\frac{q^{k(3k-1)/4}}{\bigl(q^{1/2};q^{1/2}\bigr)_k\bigl(q^{1/2};q\bigr)_k}\\
\label{rac}&=\bigl(q^{1/2};q\bigr)_\infty\sum_{k=0}^\infty\frac{q^{k(3k+1)/4}}{\bigl(q^{1/2};q^{1/2}\bigr)_k\bigl(q^{1/2};q\bigr)_{k+1}}\\
\label{rad}H(q)&=(-q;q)_\infty\sum_{k=0}^\infty\frac{(-1)^kq^{k(3k+2)}}{(q^4;q^4)_k(-q;q^2)_{k+1}}\\
\label{rae}&=\bigl(q^{1/2};q\bigr)_\infty\sum_{k=0}^\infty\frac{q^{3k(k+1)/4}}{\bigl(q^{1/2};q^{1/2}\bigr)_k\bigl(q^{1/2};q\bigr)_{k+1}}\\
\label{raf}&=\bigl(q^{1/2};q\bigr)_\infty\sum_{k=0}^\infty\frac{q^{k(3k+5)/4}}{\bigl(q^{1/2};q^{1/2}\bigr)_k\bigl(q^{1/2};q\bigr)_{k+2}}.
\end{align}
 The summations \eqref{raa}--\eqref{rae}
 are, respectively,
19, 46, 62, 97 and 44 on Slater's list (19 is equivalent to 100, and 44 to 63).
Most of them, \eqref{raa}--\eqref{rac} and \eqref{rae}, were obtained by Rogers~\cite{r1,r3}, whereas the earliest reference for \eqref{rad} is Ramanujan's lost notebook \cite[equation~(11.2.7)]{ab}.
We have not found \eqref{raf} in the literature, but as explained above it follows easily from \eqref{rae}.

The reason that Roger skips \eqref{rad} is probably that it follows by contiguity from an identity that he did know (15 and 95 on Slater's list),
\begin{equation}\label{rh}
H(q)=(-q;q)_\infty\sum_{k=0}^\infty\frac{(-1)^kq^{k(3k-2)}}{\bigl(q^4;q^4\bigr)_k\bigl(-q;q^2\bigr)_{k}},
\end{equation}
\end{subequations}
see \cite[p.\ 253]{ab} for the derivation of \eqref{rad} from \eqref{rh}.
 Not surprisingly, one can obtain \eqref{rh} by a slight modification of our proof of \eqref{rad}.
 We provide some details since it is another nice example of Andrews' constant term method.
 We first write
 \begin{equation*}
 \frac{(-1)^kq^{3k(k-2)}}{\bigl(q^4;q^4\bigr)_k\bigl(-q;q^2\bigr)_{k}}=\frac{(-1)^kq^{k(k-1)}}{\bigl(q^2;q^2\bigr)_k}\cdot \frac{(-1)^{2k}q^{2k(2k-1)/2}}{(-q;q)_{2k}},
 \end{equation*}
which leads to a constant term representation of the right-hand side of \eqref{rh},
 \begin{gather*}(-q;q)_\infty\ct \sum_{k=0}^\infty\frac{(-1)^kq^{k(k-1)}z^{2k}}{\bigl(q^2;q^2\bigr)_k}\sum_{l=-\infty}^\infty\frac{(-1)^{l}q^{l(l-1)/2}}{(-q;q)_{l}z^l}\\
 \qquad=\ct \bigl(z^2;q^2\bigr)_\infty\frac{\theta(1/z;q)}{(-qz;q)_\infty}=\ct(1+z)(z;q)_\infty\theta(1/z;q).
 \end{gather*}
 In order to apply \eqref{rct}, we use the quasi-periodicity of the theta function in the form $z\theta(1/z;q)=-\theta(q/z;q)$.
 Since the constant term does not change under $z\mapsto qz$, we can write
 \begin{align*}
 \ct(1+z)(z;q)_\infty\theta(1/z;q)&=\ct (z;q)_\infty(\theta(1/z;q)-\theta(q/z;q))\\
 &=\ct((z;q)_\infty-(qz;q)_\infty)\theta(1/z;q)=R\bigl(q^{-1};q\bigr)-R(1;q).
 \end{align*}
 This can in turn be written
 \begin{equation*}
 \sum_{k=0}^\infty\frac{q^{k^2}\bigl(q^{-k}-1\bigr)}{(q;q)_k}=\sum_{k=1}^\infty\frac{q^{k^2-k}}{(q;q)_{k-1}}=\sum_{k=0}^\infty\frac{q^{k^2+k}}{(q;q)_{k}}, \end{equation*}
and we have arrived at the left-hand side of \eqref{rh}.

Constant-term proofs of the identities \eqref{raa} and \eqref{rh} were previously given by Sills~\cite{sm}.
He uses the series \eqref{ab} rather than the original Rogers--Ramanujan identities~\eqref{rr}.

\section{Rogers' proof of the quadratic transformations}\label{rpss}

Since Rogers' papers are not easy to read, we will provide a slightly modified
version of his proof of the quadratic transformations \eqref{rtf}.
Although the notation was not used by Rogers, we find it instructive to introduce the $q$-exponential functions
\begin{align*}
e_q(x)&=\frac 1{(x;q)_\infty}=\sum_{k=0}^\infty\frac{x^k}{(q;q)_k},\qquad
E_q(x)=(-x;q)_\infty=\sum_{k=0}^\infty\frac{q^{k(k-1)/2}x^k}{(q;q)_k}.
\end{align*}
Roger defines the $q$-derivative as
\begin{equation*}(\delta_x f)(x)=\frac{f(x)-f(qx)}{x}. \end{equation*}
It is easy to see that $\delta_x e_q(tx)=t e_q(tx)$,
 and hence $f(\delta_x)e_q(tx)=f(t)e_q(tx)$ for any formal power series $f$. In particular,
\begin{equation}\label{dxe}
e_q(y\delta_x)e_q(tx)=e_q(tx)e_q(ty).\end{equation}
Identifying the coefficient of $t^k$ in \eqref{dxe} gives
\begin{equation}\label{dxh}e_q(y\delta_x)\frac{x^k}{(q;q)_k}=H_k(x,y),\end{equation}
where $H_k$ are defined by the generating function
\begin{equation}\label{eeh}e_q(tx)e_q(ty)= \frac 1{(tx,ty;q)_\infty}=\sum_{k=0}^\infty H_k(x,y)t^k.\end{equation}
In \cite{r1}, Rogers proves \eqref{dxh} in a more direct way and uses it to derive \eqref{dxe}.

The second $q$-exponential function satisfies $\delta_x E_q(tx)=t E_q(qtx)$,
which gives $\delta_x^k E_q(tx)=q^{k(k-1)/2}t^kE_q\bigl(q^ktx\bigr)$ and
\begin{equation*} e_q(y\delta_x)E_q(tx)=\sum_{k=0}^\infty\frac{q^{k(k-1)/2}y^kt^k}{(q;q)_k}\,E_q(q^ktx)
=(-tx;q)_\infty\sum_{k=0}^\infty\frac{q^{k(k-1)/2}y^kt^k}{(-tx,q;q)_k}.\end{equation*}
On the other hand, \eqref{dxh} gives
\begin{equation*}e_q(y\delta_x)E_q(tx)=\sum_{k=0}^\infty t^kq^{k(k-1)/2} H_k(x,y). \end{equation*}
Specializing $t=q$, we obtain \cite[Section~6, equation~(1)]{r2} in the form
\begin{equation}\label{rhi}\sum_{k=0}^\infty q^{k(k+1)/2} H_k(x,y)=(-qx;q)_\infty\sum_{k=0}^\infty\frac{q^{k(k+1)/2}y^k}{(-qx,q;q)_k}.
\end{equation}
By symmetry, one can interchange $x$ and $y$ on the right-hand side.

Roger now observes that the case $y=xq^{1/2}$ of \eqref{eeh} is
\begin{equation*}
 \sum_{k=0}^\infty H_k\bigl(x,xq^{1/2}\bigr)t^k=\frac 1{\bigl(tx;q^{1/2}\bigr)_\infty},\end{equation*}
which implies
\begin{equation*}
H_k\bigl(x,xq^{1/2}\bigr)=\frac{ x^k}{\bigl(q^{1/2};q^{1/2}\bigr)_k}.\end{equation*}
Using this in \eqref{rhi} and replacing $q$ by $q^2$ gives
\begin{subequations}\label{qt}
\begin{gather}
\sum_{k=0}^\infty \frac {q^{k(k+1)}x^k}{(q;q)_k}=\bigl(-q^2x;q^2\bigr)_\infty\sum_{k=0}^\infty\frac{q^{k(k+2)}x^k}{\bigl(-q^{2}x,q^2;q^2\bigr)_k}.
\end{gather}
Interchanging the role of $x$ and $y$ gives the alternative expression
\begin{gather}
\sum_{k=0}^\infty \frac {q^{k(k+1)}x^k}{(q;q)_k}=\bigl(-q^3x;q^2\bigr)_\infty\sum_{k=0}^\infty\frac{q^{k(k+1)}x^k}{\bigl(-q^{3}x,q^2;q^2\bigr)_k}.
\end{gather}
Finally, the case $y=-x$ of \eqref{eeh} is
\begin{equation*}
\sum_{k=0}^\infty H_k(x,-x)t^k=\frac 1{\bigl(t^2x^2;q^{2}\bigr)_\infty},
\end{equation*}
which gives
\begin{equation*}
H_{2k}(x,-x)=\frac{x^{2k}}{\bigl(q^2;q^2\bigr)_k},\qquad H_{2k+1}(x,x)=0.
 \end{equation*}
Using this in \eqref{rhi} gives
\begin{equation}
\sum_{k=0}^\infty \frac{q^{2k(2k+1)}x^{2k}}{\bigl(q^2;q^2\bigr)_k}=(-qx;q)_\infty\sum_{k=0}^\infty\frac{(-1)^kq^{k(k+1)/2}x^k}{(-qx,q;q)_k}.
\end{equation}
\end{subequations}
The identities \eqref{qt} are indeed equivalent to \eqref{rtf}.

\subsection*{Acknowledgements}
It is a pleasure to dedicate this paper to Stephen Milne, whose work has been one of my main sources of inspiration.
I would like to thank Andrew Sills for the reference to his master's thesis and for making it available~\cite{sm}.
I also thank the referees for several useful comments.
This research is supported by the Swedish Science Research Council, project no.\ 2020-04221.


\pdfbookmark[1]{References}{ref}
\LastPageEnding
\end{document}